\documentclass[a4paper,11pt]{article}

\usepackage{geometry,enumerate,amsmath,amssymb,latexsym,epic,theorem}
\usepackage[all]{xy}

{\theorembodyfont{\rmfamily}
\newtheorem{example}{Example}[section]}
{\theorembodyfont{\rmfamily}
\newtheorem{Def}[example]{Definition}}
{\theorembodyfont{\rmfamily} }
{\theorembodyfont{\rmfamily} }
{\theorembodyfont{\rmfamily} \newtheorem{qu}[example]{Question}}

\newtheorem{thm}[example]{Theorem}
\newtheorem{lemma}[example]{Lemma}

\newcommand{\R}{{\mathbb R}}
\newcommand{\Z}{{\mathbb Z}}

 {\theorembodyfont{\rmfamily}
\newtheorem{defn}[example]{Definition}}

\def\geq{\geqslant}

\def\ge{\geqslant}

\def\le{\leqslant}
\def\io{^{-1}}

\def\D{\mathcal{D}}
\def\s{\sigma}
\def\1{\mathbf{1}}

\def\s{\alpha}
\def\t{\beta}
\def\QED{\hfill $\Box$}
\def\L{\mathcal{L}}

\begin{document}
\title{Holonomy and monodromy groupoids\thanks{MSC2000: 18F20, 22A22, 58H05}}

\author{
 Ronald Brown  \\ School of Informatics \\ Mathematics
Division
\\ University of Wales  \\ Bangor, Gwynedd \\ LL57 1UT, U.K.
\\r.brown@bangor.ac.uk \\ \and
\.{I}lhan \.{I}\c{c}en    \\  University of  \.{I}n\"{o}n\"{u} \\
Faculty of Science and Art
\\ Department of Mathematics
\\ Malatya/ Turkey \\ iicen@inonu.edu.tr
\\ \and Osman Mucuk \\ University of Erciyes\\
Faculty of Science and Art \\Department of
Mathematics\\Kayseri/Turkey\\ mucuk@erciyes.edu.tr}

\maketitle

\begin{abstract}

We outline the construction of the holonomy groupoid of a locally
Lie groupoid and the monodromy groupoid of a Lie groupoid. These
specialise  to the well known holonomy and  monodromy groupoids of
a  foliation, when the groupoid is just an equivalence relation.
\end{abstract}

\section*{Introduction}
The holonomy and monodromy groupoids of foliations are well
known, and with their smooth structure are usually attributed to
Winkelnkemper \cite{W} and Phillips \cite{Ph}. The purpose of
this paper is to advertise the fact, due to Pradines in 1966
\cite{P1}, that these constructions are special cases of
constructions which apply to wide classes of structured
groupoids, where the foliation case is essentially that where the
groupoid is the equivalence relation determined by the leaves of
the foliation. In the final section, we suggest a number of wider
questions and possible directions for investigation, in particular
the possible relation with generalised Galois theory, and the
potentiality of higher dimensional analogues.

An important feature of Pradines' work is that these constructions
of holonomy and monodromy groupoids come with universal
properties of a local-to-global form. The association of
monodromy with a universal principle is classical, see for
example Chevalley \cite{Ch}. The monodromy principle asserts
roughly that, in a simply connected situation, for example a
simply connected group, or an equivalence relation on a simply
connected space,  a local morphism extends to a global morphism.
More generally, a local morphism can be lifted to a global
morphism on the universal cover.

The association of holonomy with a universal principle is less
well known. It is stated in terms of an adjoint pair of functors,
but not explained in detail, in \cite{P1}. It involves the notion
of what Pradines called `un morceau d'un groupo\"{\i}de
diff\'erentiables' and which we prefer to call a `locally Lie
groupoid'. This is a groupoid $G$ and a subset $W$ of $G$
containing the identities and such that $W$ has the structure of
a manifold. Conditions are imposed so that the groupoid structure
is as `smooth as possible' on $W$. There is a kind of `holonomy
principle' that, in the `locally sectionable' case (see below),
the manifold structure on $W$ extends to a Lie groupoid structure
not on $G$ but on an overgroupoid $Hol(G,W)$ of $G$, and in which
$W$ is an open subspace.

The case when $Hol(G,W)=G$ is also of interest, since this gives
a condition for the pair $(G,W)$ to be {\it extendible}. This is
used crucially to obtain a Lie structure on the monodromy
groupoid of a Lie groupoid. Thus whereas usually the holonomy
groupoid is constructed as a quotient of the monodromy groupoid,
here we regard the holonomy construction as fundamental. This
difference of approach seems of interest.

Another question arising from this work is the applicability of
the notion of locally Lie groupoid for encapsulating ideas of
local structures. It is proven by Brown and Mucuk in
\cite{Br-Muc1} that the charts of a foliation on a paracompact
manifold gives rise to a locally Lie groupoid. This process is
generalised by Brown and \.I\c{c}en in \cite{BI} to the case of a
local subgroupoid. We also note recent work of Claire Debord
\cite{De} which studies the case of singular foliations, and has
constructions whose relation to those given here would be
interesting to determine.

One aim for Pradines of this notion of what we call a locally Lie
groupoid was as a half way house between a Lie algebroid and a Lie
groupoid. We have not found  a clear statement of which  Lie
algebroids give rise to a locally Lie groupoid, but the two steps
of holonomy and monodromy groupoid were designed to model two of
the three steps in getting an essentially unique Lie group from a
Lie algebra, namely:  produce from the Lie algebra a locally Lie
group; from this produce a Lie group; finally, take the universal
cover of this Lie group. It is remarkable that Pradines'
intuitions on these steps was so strong.

The main ideas of the results and proofs for the holonomy and
monodromy groupoids were described by Pradines to Brown in the
early 1980s, and an incomplete account was written in \cite{B1}. A
full account of the holonomy construction and related material was
given in Aof's Bangor thesis \cite{A} and published in \cite{A-B}.
A full account of the monodromy construction was given in Mucuk's
Bangor thesis \cite{Mu} and published in \cite{Br-Muc1}. It should
be emphasised that this gives useful conditions for the groupoid
$M(G)$,  obtained from a Lie groupoid by taking the universal
covers of the stars of $G$ at the identities, to be given the
structure of Lie groupoid so that the projection $M(G) \to G$
gives the universal covering map on each star.

A key aspect of the construction is that $M(G)$ is initially
defined by a universal construction which ensures that it comes
with a monodromy principle on the extendibility to $M(G)$ of
certain local morphisms on $G$. The problem is to get a topology
on $M(G)$ and this, remarkably, is solved by the holonomy
construction, but in the case where the holonomy is trivial. This
seems a roundabout method. The point, however, seems to be that
the construction of the topology involves local smooth admissible
sections, and the proof that this method works seems to be no
simpler in the case of trivial holonomy than in the general case.
Thus it is important to be clear about the general method.

The use of local admissible sections for these constructions
seems essential. To see this we contrast with the group case. If
$G$ is a topological group, then left multiplication $L_g$ by an
element $g$ of $G$ maps open sets of $G$ to open sets, and in
fact $L_g$ is a homeomorphism of $G$. This is no longer the case
if $G$ is topological groupoid, for obvious domain reasons.

To remedy this situation, Ehresmann introduced the notion of
`smooth local continuous admissible section' $\sigma$ of a Lie
groupoid $G$. This is a smooth section of the source map $\alpha$
defined on some  open set $U$ of the object space $O_G$ and such
that $\beta \sigma$ maps $U$ diffeomorphically to an open set of
$O_G$. Then left multiplication $L_\sigma$ can be defined on $G$
and does map open sets of $G$ to open sets of $G$. We say that
left multiplication by an element has to be `localised', that is
`spread' to a local area. Intuitively, we regard $\sigma$ and its
associated $L_\sigma$ as a `local procedure' on the Lie groupoid
$G$.

In the case of a locally Lie groupoid $(G,W)$ there is a new
twist. We can say that $\sigma$ is smooth only if the image of
$\sigma$ lies in $W$, since only $W$ has a manifold structure. We
call such a $\sigma$ a `local procedure'. The composition in the
groupoid $G$ extends to a composition of local admissible
sections, and so such a composition can be regarded as a
`composite of local procedures', but such a composition may not
have values in $W$ and so is not a `local procedure'. In fact in
the literature, more so in physics than in mathematics, the
notion of holonomy is regarded as an iteration of local
procedures which returns to the starting point but not to the
starting value. We will see this interpreted as a germ
$[\sigma]_x$ of such a composite for which $\sigma (x)=1_x$ but
there is no neighbourhood $U$ of $x$ for which $\sigma(U)$ is
contained in $W$ and $\sigma|U$ is smooth. That is, the iteration
does not even locally give a local procedure.

The convenient formal description of the above is in terms of
inverse monoids and groupoids of germs. The nice point is that the
formal description does exactly encapsulate the intuition, and it
is the intention of this paper to convey this point.

Now we give some precise definitions.

\section{Definitions}

We fix our notation.  A groupoid consists of a set $G$ and two
functions, the source and  target maps, $\s,\t : G \to  G$  such
that $\s\t = \t , \t\s = \s$ (whence $\s^2 = \s , \t^2 = \t$, and
$\s$ and $\t$ have  the same image). We often write   $g  : \s g
\to \t g $. Further, there is a multiplication written, say, $gh
$,  for $g,h
 \in  G $, with the property that  $gh$  is defined if and only if  $\t  g = \s h$,
and then  $ \s(gh) = \s g , \t(gh) = \t h $.  The set  $\s G$   is
 called  the  set of identities, or objects, of the
groupoid  $G$, and is written $O_G$.  If  $x \in  \s G$  one
often writes $1_x$   for  $x$  to emphasise that such an  $x$
acts as an identity.  We also require associativity of the
multiplication, and the existence of an inverse to every element
of $G$.  It is often convenient to think of $O_G$ as disjoint from
$G$. Thus a groupoid is also a small category in which every
morphism is an isomorphism.

A groupoid in which $\s = \t $ is called a {\it bundle of groups},
while a groupoid in which the {\it anchor map} $(\s,\t): G \to
O_G \times O_G$ is injective is just an equivalence relation.

In order to cover both the topological and differentiable cases,
we use  the  term  ${\cal C}^r$  manifold for  $r\geq -1 $, where
the case  $r=-1$  deals with the case  of topological spaces and
continuous maps, with no  local  assumptions,  while the case
$r\geq 0$  deals as usual with  ${\cal C}^r$ manifolds and
${\cal C}^r$   maps.  Of course, a ${\cal C}^0$ map is just a
continuous map.  We then abbreviate ${\cal C}^r $ to {\em
smooth}. The terms {\em Lie group} or {\em Lie groupoid} will
then involve smoothness in this extended sense.

The following definition is due to Ehresmann~\cite{E}.

\begin{defn}  Let  $G$  be a groupoid and let  $X = O_{G}$
be a smooth  manifold. An {\em admissible local section} of $G$ is
a function $\sigma : U \rightarrow G$  from an open set in $X$
such that \begin{enumerate}[(i)] \item $\alpha \sigma(x) = x$ for
all $x \in U$; \item  $\beta \sigma(U)$ is open in $X$, and
\item  $\beta \sigma$ maps $U$  diffeomorphically to $\beta \sigma(U)$.
\end{enumerate} \label{locsecdef}
\QED
\end{defn}

Let $W$ be a subset of $G$ and let $W$ have the structure of a
smooth manifold such that $X$ is a submanifold. We say that
$(\alpha ,\beta ,W)$  is {\em  locally sectionable} if for each $w
\in W$ there is an admissible  local section $\sigma : U
\rightarrow G$ of $G$ such that (i) $\sigma\alpha (w) = w$,  (ii)
$\sigma(U) \subseteq W$ and (iii) $\sigma$ is smooth as a function
from $U$  to $W$. Such a $\sigma$ is called a {\it smooth
admissible local section.}

The following definition is due to Pradines~\cite{P1} under the
name  ``{\it morceau de groupoide diff\'erentiables}''. Recall
that if  $G$  is a groupoid then the difference map $\delta $ is
$\delta :  G \times _{\t} G \rightarrow G , (g,h) \mapstochar
\rightarrow gh^{-1} $.  \begin{Def}  A {\it locally Lie groupoid}
is a  pair $(G,W)$ consisting of a groupoid $G$ and a smooth
manifold  $W$ such that:

\noindent $G_1)$  \ \  $O_G \subseteq W \subseteq G $;

\noindent $G_2)$ \ \ $W = W^{-1} $;

\noindent $G_3)$ \ \  the set $W(\delta ) = (W \times _{\t} W)
\cap \delta ^{-1}(W)$  is open in $W \times _{\t} W$ and the
restriction of $\delta $ to  $W(\delta )$ is smooth;

\noindent $G_4)$ \  \ the restrictions to $W$ of the source and
target maps $\alpha $  and $\beta $ are smooth and the triple
$(\alpha ,\beta ,W)$ is  locally sectionable;

\noindent $G_5)$ \ \ $ W$ generates $G$ as a groupoid. \QED
\end{Def}

Note that, in this definition, $G$ is a groupoid but does not
need  to have a topology. The locally Lie groupoid $(G,W)$ is
said  to be {\it extendible} if there can be found a topology on
$G$ making it  a Lie groupoid and for which $W$ is an open
submanifold.

The main result of \cite{Br-Muc1} (which was known to Pradines)
is that a foliation $\mathcal{F}$ on a paracompact manifold $M$
gives rise to a locally Lie groupoid  $(G,W)$ where $G$ is the
equivalence relation of the leaves of the foliation, and $W$ is
constructed from a refinement of the local charts of the
foliation. In general such $(G,W)$ are not extendible. A standard
example is the foliation of the M\"obius Band $M$ by circles. In
this case the equivalence relation determined by the leaves is
not a submanifold of $M \times M$ \cite{Br-Muc1}. Foliations have
also been shown to lead to {\it local equivalence relations}
\cite{Ro1}.

Here is  an example of non extendibility due to Pradines \cite{A}.

\begin{example} \label{hol1} Consider the bundle of groups  $F$  given by the first
projection $p : \R  \times  \R   \to  \R  $, where the real line
$\R$   is considered  as a topological abelian group under
addition.  We regard  $F$  as  a  groupoid, and in fact as a
topological groupoid in the obvious sense.  Now let  $N$  be the
subgroupoid of  $F$  generated by $(x,0)$ if $x<0$ and $(x,1)$ if
$x \ge 0$.

\begin{center}
\setlength{\unitlength}{0.00053333in}
\begingroup\makeatletter\ifx\SetFigFont\undefined
\def\x#1#2#3#4#5#6#7\relax{\def\x{#1#2#3#4#5#6}}%
\expandafter\x\fmtname xxxxxx\relax \def\y{splain}%
\ifx\x\y   
\gdef\SetFigFont#1#2#3{%
  \ifnum #1<17\tiny\else \ifnum #1<20\small\else
  \ifnum #1<24\normalsize\else \ifnum #1<29\large\else
  \ifnum #1<34\Large\else \ifnum #1<41\LARGE\else
     \huge\fi\fi\fi\fi\fi\fi
  \csname #3\endcsname}%
\else \gdef\SetFigFont#1#2#3{\begingroup
  \count@#1\relax \ifnum 25<\count@\count@25\fi
  \def\x{\endgroup\@setsize\SetFigFont{#2pt}}%
  \expandafter\x
    \csname \romannumeral\the\count@ pt\expandafter\endcsname
    \csname @\romannumeral\the\count@ pt\endcsname
  \csname #3\endcsname}%
\fi \fi\endgroup
\begin{picture}(5231,4859)(0,-10)
\thicklines
 \drawline(22,2272)(4522,2272)
\drawline(2272,4822)(2272,22)
 \thinlines
\dashline{60.000}(22,2572)(4522,2572)
\dashline{60.000}(22,1972)(4522,1972)
\drawline(2272,3172)(4522,3172)
 \drawline(2272,1372)(4522,1372)
\drawline(2272,4072)(4522,4072)
\put(1850,1250){\makebox(0,0)[lb]{\smash{{{\SetFigFont{12}{14.4}{rm}$-1$}}}}}
\drawline(2272,472)(4522,472)
\put(4672,2100){\makebox(0,0)[lb]{\smash{{{\SetFigFont{41}{49.2}{rm}$\}$}}}}}
\put(5047,2197){\makebox(0,0)[lb]{\smash{{{\SetFigFont{12}{14.4}{rm}$W'$}}}}}
\put(2047,3060){\makebox(0,0)[lb]{\smash{{{\SetFigFont{12}{14.4}{rm}$1$}}}}}
\put(2047,2030){\makebox(0,0)[lb]{\smash{{{\SetFigFont{12}{14.4}{rm}$0$}}}}}
\end{picture}

Figure 1

\end{center}

Let  $G$   be the quotient groupoid  $F/N $, and let  $q : F  \to
G$  be  the quotient morphism of groupoids.  Then the stars
$\alpha \io (x)$ of $G$ are bijective with $\R$ if $x <0$ and
with $\R/\Z$ if $ x \ge 0$. Let $W'$ be the subspace $\R \times
(-1/4, 1/4 )$ of $F$, and let  $W = qW'$. The topology on $W'$
may easily  be transferred to  a topology on $W$   so that $(G,W)$
becomes a locally Lie groupoid. However, it is not possible to
extend this topology so as to get even a topological groupoid
structure on  $G$, for which $W$   is an open  subspace. This can
be seen by noting that the section $s$ of the map $\alpha$ of $G$
given by $x \mapsto q(x,1/8)$  is continuous but $9s$ is not.
Instead, there is another groupoid $H = Hol(G,W)$, called the
holonomy groupoid of the locally topological groupoid $(G,W)$,
which is a topological groupoid, and which contains $W$   as an
open subspace.  The groupoid $H$ is equipped with a surjective
morphism   $\phi  : H \to G$ which is the identity on objects. In
this case the kernel of $\phi$ is non trivial only at $0$ and is
there of the form $\{ 0\}  \times \Z$. \QED
\end{example}

\begin{example} There is a variant of this last example  in  which $F$   is as
above, but this time  $N$  is the union of the groups  $\{ x\}
\times   (1  + |x|)\Z$   for all  $x  \in  \R$.  One defines
$W$    as  before,  but  this time considers  $W$   as a
differential  manifold.   The  topological structure on $W$ can
be  extended  to  give  a  topological  groupoid structure on the
quotient  $G = F/N$.  The differential structure, however, cannot
be so extended, because the section $s$ given as in the previous
example is such that $9s$ is not smooth. In this case one gets a
differential holonomy groupoid, with a projection morphism $\phi
: H  \to G$ whose kernel is as in the previous example. One can
get similar examples with varying degrees of differentiability
considered. In this and the previous example, the  holonomy
groupoids constructed are non-Hausdorff topological  (or Lie)
groupoids. \QED
\end{example}

\section{The holonomy construction}
The main result of Aof and Brown \cite{A-B} is a version of
Th\'{e}or\`{e}me 1 of  Pradines~\cite{P1} and was stated in the
topological case.  In the smooth case it states:

\begin{thm} {\em (Pradines~\cite{P1}, Aof and
Brown~\cite{A-B}) (Globalisability theorem)}  Let $(G,W)$  be a
locally Lie  groupoid. Then there is a Lie groupoid $H$,  a
morphism $\phi : H \rightarrow G$   of groupoids and an embedding
$i : W \rightarrow H$ {\it of} $W$ to an open neighborhood of
$O_{H}$ such that the following conditions are satisfied.

\noindent {\em i)} $\phi $ is the identity on objects, $\phi i =
id_{W} ,  \phi ^{-1}(W)$ is open in $H$, and the restriction
$\phi _{W} : \phi^{-1}(W)  \rightarrow W$  of $\phi $  is smooth;

\noindent {\em ii) } if A is a Lie groupoid and $\xi :  A
\rightarrow G$  is a morphism of groupoids such that:

a) $\xi $ is the identity on objects;

b) the restriction $\xi _{W} : \xi ^{-1}(W) \rightarrow W$  of
$\xi $  is smooth and $\xi ^{-1}(W)$ is open in A and generates A;

$c)$  the triple $(\alpha _{A} , \beta _{A} , A)$  is locally
sectionable,

\noindent then there is a unique morphism $\xi ^\prime : A
\rightarrow H$   of Lie  groupoids such that $\phi \xi ^\prime =
\xi $   and $\xi ^\prime a = i\xi a$ for $a \in \xi ^{-1}(W)$.
\end{thm}

The groupoid $H$ is called the {\it holonomy groupoid} $Hol(G,W)$
of the  locally Lie groupoid $(G,W)$. It is thus the minimal
overgroupoid of $G$ which can be made into a Lie groupoid with
$W$ as open subspace.

We should also say that Pradines actually states more since his
is a theorem on germs of such $(G,W)$. So there is still more
work to be done on giving a full account of this result and
illustrating it with examples.

\noindent {\bf Sketch of the proof of Theorem 2.1}

An important construction due to Ehresmann is a multiplication on
the set   $\Gamma (G)$ of local admissible sections of $G$ in
which if  $x  \in  X$
   $$         ( \sigma \tau )(x) = ( \sigma x)(\tau \t \sigma x).$$
With this multiplication,  $\Gamma (G)$  is a monoid, and in fact
an  inverse monoid, in the sense that every   $\sigma$ has  a
unique (generalised) inverse  $\sigma '$ such that $$\sigma
\sigma ' \sigma  =  \sigma  , \qquad  \sigma ' \sigma  \sigma ' =
\sigma ' . $$ Since  $\sigma 'x  =  (  \sigma  (\t  \sigma )\io
x)\io $, we write  $\sigma \io$   for   $\sigma ' $.  An
important reason for introducing these sections is that if $G$
is  a topological groupoid, then translation by a continuous
local admissible section  does  map  open sets of $G$   to open
sets of  $G$.

Let  $\Gamma ^c (W)$  be the subset of  $\Gamma (G)$ consisting of
local admissible sections  which  (i) have  values in   $W$ and
(ii) are smooth.  Of course the first condition  is necessary for
the second condition to make sense.  Let $\Gamma ^c (G,W)$ be the
sub-inverse monoid of $\Gamma (G)$ generated by $\Gamma ^c (W) $.
At this stage it is convenient to assume that $W = W\io $. It is
proved in \cite{A} that this is no loss of generality.

Now let  $J^c (G,W)$  be the sheaf of germs of  the  elements  of
$\Gamma ^c (G,W)$, and let  $J^c (W)$  be the sheaf  of  germs  of
the elements  of   $\Gamma ^c (W)$. The germ of a local section
$\sigma$   at the point  $x$  of  its  domain  is written $[
\sigma ]_x  $.  Then  the  inverse monoid  structure  on $\Gamma
(G)$ induces on  $J^c (G,W)$ the structure of groupoid, in which
                 $$[ \sigma ]_x [\tau ]_y = [ \sigma \tau ]_x$$
is defined if and only if  $y = \t \sigma x $.

The sets  $\Gamma ^c (W)$  and  $J^c (W)$  have a r\^{o}le as
codifying  a {\it local procedure}.   The  inverse  monoid $\Gamma
^c  (G,W)$ and  the  groupoid $J^c (G,W)$ then codify the {\it
iteration of local procedures}.  It is in this sense that we are
dealing with local-to-global techniques.

There is a morphism of groupoids, the `final map',    $\psi  : J^c
(G,W) \to G , [ \sigma ]_x  \mapsto  \sigma x $, which is the
identity on objects.  We set
        $$            J_0 = J^c (W)  \cap  (Ker  \psi ) ,
$$ so that  $J_0$  consists  of  germs   $[  \sigma  ]_x$    of
continuous  local admissible sections   $\sigma$   with values in
$W$   and such that   $\sigma x = 1_x  $. The aim is to define the
holonomy groupoid of the locally Lie groupoid  $(G,W)$ to be the
quotient groupoid
                           $$Hol(G,W) = J^c (G,W)/J_0 .$$
For this we need to prove:

\begin{lemma}   The set  $J_0$  is a normal subgroupoid of  $J^c (G,W) $.
\end{lemma}

The main point of the proof is that because of the definitions of
$J_0$ and of  $J^c (G,W)$  one has only to check that if  $[\rho
]_x ,[ \sigma ]_x   \in J_0$ and $[\tau ]_x  \in  J^c (W) $, then
$[\rho \sigma \io ]_x   \in  J_0 $, and $[\tau ]_x [ \sigma  ]_x
[\tau ]_x  \io    \in   J_0  $.  This  follows from continuity
considerations and the facts that
            $$(\rho  \sigma )x = 1_x  = (\tau  \sigma \tau \io )x
.$$

Let  $p : J^c (G,W)  \to  Hol(G,W)$  be the  quotient  morphism.
We  write $H$   for  $Hol(G,W)$  and write  $\langle  \sigma
\rangle _x$ for   $p([  \sigma ]_x )$. Thus  $H$   is a groupoid.
Note that the morphism   $\psi  : J^c (G,W)   \to   G$ induces a
morphism which we write   $\phi  : H  \to   G  ,  \langle \sigma
\rangle _x \mapsto \sigma x $.  For this morphism to be
surjective, it is sufficient to  assume that  $W$   generates
$G$   as a groupoid, and that for every  element   $w$ of  $W$
there is a continuous admissible local section of  $G$    through
$w$.

Let  $f  \in  \Gamma ^c (G,W)$.  We define a partial function
$\chi  _f : W \to  H$, by
  $$w \mapsto  \langle f\rangle _x \langle  \sigma_w \rangle _x  ,$$
where   $\sigma_w$   is an admissible local section of  $s$
through $w$.   Again, one has to assume that such a section
exists for all  $w  \in  W $, and  one  has to prove that this
value  is independent  of  the  choice  of  local  section
$\sigma_w  $, and that   $\chi _f$  is injective with domain an
open  subset  of $W$.

A key lemma is that if $f,g  \in  \Gamma ^c (G,W)$ then  $( \chi
_f )\io ( \chi _g )=L_h$, left multiplication by the section $f\io
g$. This shows that $( \chi _f )\io ( \chi _g )$ maps an open set
of $W$ diffeomorphically to an open set of $W$. This algebraic
format for the change of charts is also convenient for proving
$Hol(G,W)$ becomes a Lie groupoid, see \cite{A-B}.

We also need that every element of the holonomy groupoid arises
in this way, and for this we also need that  $W$ generates  $G$.
Such an assumption is in practice not so great a restriction.  A
result of Pradines (compare \cite[Proposition 1.5.16]{A}) is that
if each star  $\s \io(x)$ of $G$   meets  $W$ in  a connected set,
then any open neighbourhood of  $X$   in $W$   generates $G$.

These results allow the   $\chi _f$   for all  $f  \in  \Gamma ^c
(G,W)$  to be used as charts for a topology on the holonomy
groupoid  $H$. Notice that every element of  $H$   is of the
form   $\chi _f (x)$ for some  $f  \in \Gamma ^c (G,W)$ and  $x
\in  \D_ f$. Consequently,  given   $f   \in   \Gamma ^c (G,W) $,
the function $x \mapsto   \chi _f x$  for $x  \in  \D_ f$  is a
continuous admissible local section of  $H$. Also,  $H$   is
generated as a groupoid by    $\chi _{\mathbf 1} (W)$ where
$\mathbf 1$  here denotes the identity section with domain  $X$.
This completes the sketch proof. \QED

Readers of the Bourbaki account for Lie groups  (\cite{Bo}
p.210)  may  be puzzled by the lack of a condition involving
conjugacy, of  the type  that for all  $g  \in  G$  there is an
open neighbourhood $U$  of   $\s g$   such  that $gUg\io$  is
contained in  $W$. Pradines argues (private communication) that
in the first  place this condition  is  unrealistic,  since  it
involves `global' elements  $g$  of  $G$.  In the second place,
this  condition is not needed, by virtue of the assumptions on
generation.

The above construction can be followed through to give the
results  of Examples 5 and 6.

There  is  a  surprising  application   of   the   holonomy
groupoid construction, namely  to  give  a  condition  that  a
locally Lie groupoid  $(G,W)$   is extendible,  i.e. determines a
topology  on   $G$ making it a Lie groupoid for which $W$ is an
open subspace.   In terms of previous notation, this condition is
simply that  $Ker  \; \psi$    is contained in $J^c (W)$, which
is equivalent to the condition that if    $\sigma$ is any product
of admissible continuous local sections about  $x$ each with
values in  $W$, and   $\sigma (x) =  1_x$,  then some
restriction  of $\sigma$   to a neighbourhood of $x$ has values
in $W$ and is smooth. It is not clear that there is any easier
proof of this extendibility result than  that  obtained from the
construction of the holonomy groupoid.

This extendibility result is  used,  as  suggested  by  Pradines
(see \cite{B1}), in constructing a topology on the monodromy
groupoid of a topological groupoid.  The basic method is as
follows.

Let now  $G$   be a topological groupoid and let  $W$   be an open
subset of $G$   containing the identities.  The groupoid structure
on  $G$   makes  $W$ into a pregroupoid, by which is meant that
the product    $uv$   of  two  elements $u,v$  of  $W$ is not
always defined (in  $W$).  There is  a  standard  way of making
any pregroupoid   $W$    into  a  groupoid   $M$    with a
morphism of pregroupoids  $i : W  \to  M$  such that any
pregroupoid morphism from  $W$   to a groupoid  $K$   extends
uniquely  to  a morphism   $M   \to   K  $.   Since $W$ embeds in
a groupoid (namely  $G$), the morphism  $i :  W   \to   M$   is an
embedding.  Methods of \cite{D-L} may be extended to show  that
under suitable local conditions on  $G$, the topology on $W$ may
be extended to a topology on each  $s_M\io x , x  \in  X $, such
that each projection   $s_M\io x   \to s_G\io x$ is a universal
cover. The  previously  mentioned  universal property now gives a
version of the classical Monodromy Principle \cite{Ch}, but stated
in terms of groupoids, rather than equivalence relations or groups
as in \cite{Ch}.

The problem is now to make  $M$   into a topological  groupoid so
that the universal property yields  a  continuous  morphism  on
$M$ if $K$   is  a topological groupoid and the pregroupoid
morphism $W \to  K$  is continuous. The surprising, but simple to
prove,  result is that the pair $(M,W)$ satisfies the condition
for extendibility  stated above, basically because $G$ is already
a Lie groupoid. So the holonomy method outlined above is used to
extend the topology on $W$   to a topology on $M$, assuming that
$G$ has enough continuous admissible local sections.  (This is
Pradines' method for Th\'eor\`eme  2 of \cite{P1}, 1966, explained
to the Brown in 1981.)

The monodromy groupoid construction yields the homotopy groupoid
of a foliation, discussed in  \cite{Ph}. It also yields this
groupoid with the universal property of globalising a morphism
defined locally. Once again we see a local-to-global feature
which fits naturally into  the context of groupoids.

These methods also give an answer to the following question:
given a Lie groupoid $G$, let $Cov(G)$ be the union of the
universal covers at $1_x$ of the star of $G$ at $x$, for all $x
\in G$. Let $q:Cov(G)\to G$ be the projection. Assume we can
choose a neighbourhood $W$ of $O_G$ so that the inclusion $W \to
G$ lifts to $W \to Cov(G)$. The monodromy principle then yields a
morphism of groupoids $\phi: M(G) \to Cov(G)$ which is continuous
on stars. But $p:M(G) \to G$ is a covering map on each star, and
so $\phi$ is a bijection on each star, and hence is an
isomorphism. This isomorphism induces a topology on $Cov(G)$
making it a topological groupoid, or, in appropriate
circumstances, a Lie groupoid. Such a construction is given by
Mackenzie in \cite{M} by a different method, in the locally
trivial case.

The full details of the above arguments are given in \cite{B-M}.

The use of the monodromy groupoid and $Cov(G)$ also enables us to
explain the holonomy groupoid of Example \ref{hol1}. The monodromy
groupoid of $Hol(G,W)$ is the original groupoid $F$ and so
$Hol(G,W)$ is isomorphic to the quotient of $F$ by the
subgroupoid generated by $(x,0)$ for $x \le 0$ and $(x,1)$ for $x
>0$.
\section{Local subgroupoids}

There is considerable work on local equivalence relations, part
of the motivation being that a foliation on a manifold $M$
determines a local equivalence relation on $M$ \cite{Ro1}. Now an
equivalence relation on $M$ is just a  subgroupoid of the
indiscrete groupoid $M \times M$ which has the object set $M$
(this is also known as a {\it wide} subgroupoid of $M \times M$).
It thus seems natural to consider an arbitrary groupoid $Q$ with
object set $M$ and to consider the sheaf $p:\L_Q \to M$
associated to the presheaf $U \mapsto L_Q(U)$ where $L_Q(U)$ is
the set of wide subgroupoids of $Q|U$. This notion is studied in
\cite{BI,BIM}. In \cite{BI} there are given conditions on a local
subgroupoid of a Lie groupoid so that it leads to a locally Lie
groupoid and hence to holonomy and monodromy groupoids. In
particular, this leads to a monodromy principle for local
subgroupoids.

\section{Questions}

\begin{qu}We have already mentioned the question of extending  the work
on holonomy and monodromy to germs, thus giving a complete account
of the theorems in the first of Pradines' notes \cite{P1}, which
are stated as the existence of adjoint functors. Some remarks on
this are given in \cite{P2}.
\end{qu}
\begin{qu} It would be interesting to know (i) how useful is the
notion of locally Lie groupoid in formulating local properties,
and (ii) what is its relation to the notion of Lie algebroid.
\end{qu}
\begin{qu}
The following question could be of interest.  Grothendieck has
developed extensive work on the fundamental group in the context
of  algebraic geometry.  The notion of monodromy is also often
vital in these  arithmetic questions. Can the above approach to
monodromy and covering spaces  be  of use in  these arithmetic
problems?   This  would be   an interesting further vindication
of Pradines' approach, and is related to the next question.
\end{qu}
\begin{qu} It is known that covering spaces and Galois theory are closely related,
see for example \cite{DD,BJvkt}. The last paper relates the
generalised Galois theory of Janelidze \cite{J} to covering space
theory. It would be very interesting to tie in notions of
monodromy for groupoids with these broader aspects of Galois
theory and of descent \cite{JT}.
\end{qu}
\begin{qu} For the present writers, the most intriguing, and possibly the most difficult,
question is that of higher dimensional analogues of these
results. Background to the idea that multiple groupoids form
candidates for `higher dimensional groups' is given in \cite{B4}.
A starting point was  that  since higher homotopy groups are
abelian because `group objects in groups' are just abelian groups,
it is therefore natural to look at objects of the type of `group
objects in groupoids' or `groupoid objects in groupoids'. These
are more complicated objects than groups, and the complication of
$n$-fold groupoids increases directly with $n$. Indeed, it is
known that $n$-fold groupoids model homotopy $n$-types. Such
$n$-dimensional structures lend themselves to the consideration
of `algebraic inverses to subdivision'; since subdivision is a
fundamental process in local-to-global questions, the possibility
of detailed algebraic control over the inverse process in certain
circumstances would be expected to lead to surprising new
results, and involving essentially non abelian considerations.
These objects  do arise naturally in homotopy theory, where they
lead to new algebraic constructions such as a non abelian tensor
product of groups and to calculations in homotopy theory not
possible by other means \cite{BL,B4}.  These algebraic objects,
or analogous ones, also arise in many other algebraic and
geometric situations \cite{BMac,BJ,M2,MX}.

Thus it is natural to consider the possibility of higher
dimensional forms of holonomy and monodromy. A tentative step in
this direction is given in \cite{BIdoub}, which covers part of
\cite{I}.  The basic intuition is that for a groupoid an
admissible section can also be considered as a homotopy. A
reasonable generalisation of an admissible section should
therefore be a notion of a homotopy, i.e. a deformation. This
notion exists for various forms of double groupoid. Thus the
existence of multiple geometric structures (double foliations,
foliated bundles, etc.) should in principle be properly reflected
by multiple algebraic structures.

It is surely intuitively significant in this respect that
multiple categories arise in the context of concurrency in
computer science, where the multiple processors are thought of as
each giving another time dimension. The algebraic analysis seems
naturally to  involve a generalisation of the notion of free
category on a graph to a certain notion of a free cubical
$\omega$-category on a cubical set. The analysis of this situation
is still incomplete, but is studied in \cite{Go,Ga}.

It is possible that a description of the relation between holonomy
in the sense of this paper and holonomy for principal bundles
with connection, and hence the relation with curvature, requires
some higher dimensional algebraic treatment.
\end{qu}

\section*{Acknowledgements}  We  would like to thank J.Pradines
for sharing his ideas on holonomy and monodromy groupoids, and
Mohammed Aof for his development of these ideas.

\end{document}